\documentclass{article}
\usepackage{latexsym}
\usepackage{amssymb}
\usepackage{amsmath}
\usepackage{amscd}
\usepackage[only,twlrm,ninrm,sevrm]{rawfonts}
\usepackage[latin1]{inputenc}
\numberwithin{equation}{section} \textwidth=5.45in
\newtheorem{Th}{Theorem}[section]

\newtheorem{Le}{Lemma}[section]
\newtheorem{Remark}{Remark}[section]

\newcommand{\R}{\mathrm{I\!R\!}}

\newcommand{\N}{\mathrm{I\!N\!}}

\newcommand{\var}{\varepsilon}
\newcommand{\dis}{\displaystyle}
\newcommand{\rig}{\rightarrow}
\newcommand{\cha}{\widehat}
\newcommand{\righ}{\rightharpoonup}

\newcommand{\fim}{\blacksquare}

\begin{document}
\title{The effect of the domain topology on the number of positive solutions of an elliptic Kirchhoff problem}
\author{Jo\~ao R. Santos Junior\footnote{Supported by CAPES - Brazil - 7155123/2012-9}\\
Faculdade de Matem\'atica\\
Universidade Federal do Par\'a\\
66.075-110-Bel\'em-Par\'a - Brazil\\}

\date{}
\maketitle \vspace*{-25pt}


\begin{abstract}
Using  minimax methods and Lusternik-Schnirelmann theory, we study
multiple positive solutions for the Schr\"{o}dinger - Kirchhoff
equation
$$
M\left(\dis\int_{\Omega_{\lambda}}|\nabla
u|^{2}dx+\dis\int_{\Omega_{\lambda}}u^{2}dx\right)\left[-\Delta u +
u \right]= f(u)
$$
in $\Omega_{\lambda} = \lambda\Omega$. The set $\Omega \subset
\mathbb{R}^3$ is a smooth bounded domain, $\lambda>0$ is a
parameter, $M$ is a general continuous function and $f$ is a
superlinear continuous function with subcritical growth. Our main
result relates,  for large values of $\lambda$, the number of
solutions with the least number of closed and contractible in
$\overline{\Omega}$ which cover $\overline{\Omega}$. \vspace{0.5cm}


\end{abstract}


\section{Introduction}

In this paper we study multiple positive solutions for the following
problem
$$
 \left\{
\begin{array}{rcl}
\mathcal{L}u=f(u), \ \Omega_{\lambda}\\
u>0, \ \Omega_{\lambda}\\
u=0, \ \partial\Omega_{\lambda}
\end{array}
\right.\leqno{(P_{\lambda})}
$$
where $\Omega \subset \mathbb{R}^{3}$ is a smooth bounded domain,
$\lambda>0$ is a parameter, $\Omega_{\lambda}:=\lambda\Omega$ is an
expanding domain and $\mathcal{L}$ is the nonlocal operator given by
$$
\mathcal{L}u=M\left(\dis\int_{\Omega_{\lambda}}|\nabla
u|^{2}dx+\dis\int_{\Omega_{\lambda}}u^{2}dx\right)\left[-\Delta u +
u \right].
$$

In 1883, Kirchhoff \cite{kirchhoff} established the equation
$$
\rho\displaystyle\frac{\partial^{2}u}{\partial
t^{2}}-\biggl(\displaystyle\frac{P_{0}}{h}+\displaystyle\frac{E}{2L}\int^{L}_{0}\biggl|\displaystyle\frac{\partial
u}{\partial
x}\biggl|^{2}dx\biggl)\displaystyle\frac{\partial^{2}u}{\partial
x^{2}}=0 \leqno{(K)}
$$
where  $L$ is the length of the string, $h$ is the area of
cross-section, $E$ is the Young modulus of the material, $\rho$ is
the mass density and $P_{0}$ is the initial tension. This model was
proposed to modify the classical d'Alembert's wave equation,
assuming a nonlinear dependence of the axial strain on the
deformation of the gradient.

Owing to its importance in engineering, physics and material
mechanics, a considerable effort has been devoted during the last
years to the study the generalization of the stationary equation
associated with problem $(K)$. With no hope of being thorough, we
mention some papers regarding the study of this class of problems:
\cite{alvescorreama}, \cite{jmaa}, \cite{He}, \cite{Li}, \cite{ma},
\cite{Wang 1}, \cite{Wang 2}, \cite {Wu} and reference therein. For
an excellent didactic about this class of problems we cite
\cite{Azzollini} and for an overview of non-local problems we cite
\cite{Giovany}.

Problem $(P_{\lambda})$ is a generalization of the stationary
problem associated with problem $(K)$. Before stating our main
result, we need the following hypotheses on the functions $M$ and
$f$.

The continuous function  $M: \mathbb{R}_{+} \to \mathbb{R}_{+}$  and
the nonlinearity $f: \mathbb{R} \to \mathbb{R}$ satisfies the
following conditions:

\begin{description}

\item[($M_{1}$)] There is $m_{0}>0$ such that $M(t)\geq m_{0}, \ \forall t\geq 0$.

\item[($M_{2}$)] The function $t\mapsto M(t)$ is increasing.

\item[($M_{3}$)] The function $t\mapsto\displaystyle\frac{M(t)}{t}$ is decreasing.

\end{description}

A typical example of function verifying the assumptions
$(M_{1})-(M_{3})$ is given by $\dis M(t)=m_{0}+bt$, with $m_{0}>0$
and $b>0$. More generally, each function of the form $\dis
M(t)=m_{0}+bt+\dis\sum_{i=1}^{k}b_{i}t^{\gamma_{i}}$ with $b_{i}\geq
0$ and $\gamma_{i}\in (0,1)$ for all $i\in \{1,2,\ldots, k\}$
verifies the hypotheses $(M_{1})-(M_{3})$.

Now we give an example  of a continuous but non-differentiable
function that satisfies such hypotheses. Let   $m_{0}, b_{0}, b_{1}$
and $t_{0}$ be positive constants such that $b_{0}\neq b_{1}$ and
$t_{0}<\frac{m_{0}}{b_{1}-b_{0}}$ if $b_{0}<b_{1}$. We define the
continuous function
$$
M(t)= \left\{
\begin{array}{rcl}
m_{0}+b_{0}t, \ \mbox{if} \ 0\leq t\leq t_{0}\\
m_{0}+(b_{0}- b_{1})t_{0}+b_{1}t, \ \mbox{if} \ t_{0}\leq t\\
\end{array}
\right.
$$
Since that $b_{0}\neq b_{1}$, we have that  $M$ is
non-differentiable in $t_{0}$. Using the same reasoning, we can
build  continuous functions that are not differentiable in a finite
number of points.

We assume that the locally Lipschitz continuous function $f$ vanishes in
$(-\infty, 0)$ and verifies
\begin{description}
\item[($f_{1}$)] $$\lim_{t\to 0^{+}}\frac{f(t)}{t^{3}}=0. $$

\item[($f_{2}$)] There is $q\in (4,6)$ such that

$$
\lim_{t\to \infty}\frac{f(t)}{t^{q-1}}=0.
$$

\item[($f_{3}$)]  There is $\theta\in (4,6)$ such that
$$
0<\theta F(t)\leq f(t)t, \ \forall t>0,
$$
where $F(s)=\displaystyle\int_{0}^{s}f(t)dt$.

\item[($f_{4}$)] The application $$t\mapsto\frac{f(t)}{t^{3}}$$ is nondecreasing in $(0,\infty)$.

\end{description}

A typical example of locally Lipschitz continuous function verifying the assumptions $(f_{1})- (f_{4})$ is given by 
$$
f(t)=\displaystyle\sum_{i=1}^{n}c_{i}(t^{+})^{q_{i}-1}
$$
with $c_{i}\geq 0$ not all zero and $q_{i}\in [\theta, 6)$ for all $i\in \{1,2,\ldots, n\}$. Moreover, a very simple example of 
non-defferentiable function verifying these hypotheses is given by
$$
f(t)= \left\{
\begin{array}{rcl}
c(t^{+})^{q_{2}}, \ \mbox{if} \ 0\leq t\leq 1\\
c(t^{+})^{q_{1}}, \ \mbox{if} \ t\geq 1,\\
\end{array}
\right.
$$
where $c>0$ and $4<\theta\leq q_{1}<q_{2}<6$.

The main result of this paper is:

\begin{Th}\label{Theorem 1.1}
Suppose that the function $M$ satisfies $(M_{1})-(M_{3})$ and the
function $f$ satisfies $(f_{1})-(f_{4})$. Then there exists
$\lambda^{*}>0$ such that, for each $\lambda\in[\lambda^{*}, \infty)$, the
problem $(P_{\lambda})$ has at least $cat \Omega$ positive weak solutions. 
Moreover, if $cat\Omega>1$ then $(P_{\lambda})$ has at least $cat\Omega + 1$
weak solutions.
\end{Th}

For more informations about the Lusternik - Schnirelmann category, we refer to $\cite{BC}$ and $\cite{BC 2}$.

\vspace{.3cm}

When $M=1$ we have the following problem
$$
 \left\{
\begin{array}{rcl}
-\Delta u +u =f(u), \ \Omega_{\lambda}\\
u>0, \ \Omega_{\lambda}\\
u=0, \ \partial\Omega_{\lambda}
\end{array}
\right.\leqno{(BC)}
$$
that was studied  first by Benci and Cerami in \cite{BC}. In order
to obtain multiple solutions for this problem, these authors made
comparisons between the category of some sublevel sets of the
functional associated to the problem and the category of the domain
$\Omega_{\lambda}$. After this excellent paper, appeared several
generalizations. A version of this problem for a class of
quasilinear equation can be seen in \cite{AGU}. In \cite{AFF} there
is a version considering the Schr\"{o}dinger operator in the
presence of magnetic potential. The case with p-Laplacian operator
is in \cite{Alv1}. In all these works the nonlinearity $f$ is
$C^{1}$ class, because in the arguments used was important the
regularity of the Nehari manifold associated to $(P_{\lambda})$.

Problem $(P_{\lambda})$ is a nonlocal version of the $(BC)$
considering the Kirchhoff operator. But the presence of Kirchhoff
operator with $M$ and $f$ only continuous imply that several
estimates used in \cite{BC}, \cite{Alv1}, \cite{AGU} and \cite{AFF}
cannot be repeated for the functional energy associated to
$(P_{\lambda})$. To overcome this difficult we use an argument that
can be found in \cite{Szulkin} and  \cite{Wang 1}, and we introduce some Lemmas,
as for example, Lemmas $\ref{PVE}$ and $\ref{compacityLemma}$. However, due to
the presence of the function $M$, some estimates more refined are
need, such as in the study of the limit problem and in the Lemma
\ref{lema-assintotico1}.

An important point in this type of arguments is the existence
of solution of a limit problem. In our case, the limit problem is
given by
$$
 \left\{
\begin{array}{rcl}
\mathcal{L}_{\infty}u=f(u), \ \R^{3}\\
u>0, \ \R^{3}.
\end{array}
\right.\leqno{(P_{\infty})}
$$

This result was proved in \cite{Alves and Figueiredo}. In our paper
we show this result of existence  considering a less restrictive set of assumptions
about $f$ and $M$.

The paper is organized as follows. In the section 2 we study the
existence of solution of the limit problem and we prove a compactness result of the
Nehari manifold associated to the functional of the limit problem.
This study was not necessary in \cite{Alves and Figueiredo}. In the
section 3 we study the behavior of minimax levels from the
functional associated to the problem $(P_{\lambda})$. The main
result is proved in the section 4.


\section{The limit problem}

An important result that we shall use in this work is related to the
existence of a positive ground-state solution for the problem
$$
 \left\{
\begin{array}{rcl}
\mathcal{L}_{\infty}u=f(u), \ \R^{3}\\
u>0, \ \R^{3},
\end{array}
\right.\leqno{(P_{\infty})}
$$
where
$$
\mathcal{L}_{\infty}u=M\left(\dis\int_{\R^{3}}|\nabla
u|^{2}dx+\dis\int_{\R^{3}}u^{2}dx\right)\left[-\Delta u + u \right].
$$
More precisely, we are concerned with the existence of a positive function $u \in
H^{1}(\mathbb{R}^{3})$ verifying
$$
I_{\infty}(u)= c_{\infty} \ \ \mbox{and} \ \ I'_{\infty}(u)=0,
$$
where
$$
I_{\infty}(u)=\frac{1}{2}\cha{M}(\| u\|^{2})-\dis\int_{\R^{3}}F(u) \
dx,
$$
$\widehat{M}(t)=\displaystyle\int^{t}_{0}M(s) \ ds$ and $c_{\infty}$
denotes the minimax level of the mountain pass theorem associated to
the functional $I_{\infty}$ and given by
$$
c_{\infty}=\inf_{\gamma\in \Gamma}\sup_{t\in
[0,1]}I_{\infty}(\gamma(t)),
$$
with $\Gamma=\{\gamma\in C([0,1], H^{1}(\mathbb{R}^{3}): \gamma(0)=0
\ \mbox{and} \ I_{\infty}(\gamma(1))<0\}$.

It is not difficult to check that $I_{\infty}$ is a $C^{1}$ functional,
$$
I'_{\infty}(u)v=M(\| u\|^{2})(u,v)-\dis\int_{\R^{3}}f(u)v \ dx, \
\forall u,v\in H^{1}(\R^{3})
$$
and its Nehari manifold is given by
$$
\mathcal{N}_{\infty}=\{u\in H^{1}(\R^{3})\backslash\{0\}:
I'_{\infty}(u)u=0\}.
$$
Here, $(,)$ and $\|.\|$ denote, respectively, the standard inner product of $H^{1}(\mathbb{R}^{3})$ and its induced norm. 

In \cite[Theorem 2.5]{Alves and Figueiredo} was proved that problem
$(P_{\infty})$ has a positive solution. Since that in our paper we
have a smaller number of the hypotheses, we give some details about
this result.

\subsection{Existence of positive ground-state solution for the limit
problem}

Notice from $(M_{3})$, there is a positive constant $K$ such that
\begin{eqnarray}\label{crescimentoM}
M(t)\leq K + M(1)t,
\end{eqnarray}
for all $t\geq 0$. Thus, the growth condition  (\ref{crescimentoM})
and $(M_{1})$ allow us to use \cite[Lemma 2.1,
Lemma 2.2]{Alves and Figueiredo} and to conclude that functional
$I_{\infty}$ has the Mountain Pass geometry and by a version of
Mountain Pass Theorem (see [$\cite{Willem}$, Theorem 1.15]), there
is  a $(PS)_{c_{\infty}}$ sequence for the functional $I_{\infty}$,
that is, there is a sequence $(u_{n})\subset H^{1}(\R^{3})$ such
that

$$
I_{\infty}(u_{n})\rig c_{\infty} \ \mbox{and} \
I_{\infty}'(u_{n})\rig 0.
$$

Using \cite[Lemma 2.4]{Alves and Figueiredo} we can assume that the
weak limit of a $(PS)_{c_{\infty}}$ is nontrivial. Indeed, suppose that 
$u_{n}\righ 0$. Since $u_{n}\not\rightarrow 0$, there are $(y_{n})\subset \R^{3}$ and $R,\beta>0$ such
that
$$
\liminf_{n\to \infty}\int_{B_{R}(y_{n})}u_{n}^{2} \ dx\geq \beta>0.
$$
Considering $v_{n}(x)=u_{n}(x+y_{n})$, we can prove that $(v_{n})$
is a $(PS)_{c_{\infty}}$ sequence for the functional $I_{\infty}$,
$(v_{n})$ is bounded in $H^{1}(\R^{3})$ and there is $v\in H^{1}(\R^{3})$ non-trivial,
with $v_{n}\rightharpoonup v$ em $H^{1}(\R^{3})$.

In the next Proposition we obtain a positive ground-state solution
for the autonomous problem $(P_{\infty})$.

\begin{Th}\label{Theorem 3.1}
Let $(u_{n})\subset H^{1}(\R^{3})$ be a $(PS)_{c_{\infty}}$ sequence
for $I_{\infty}$. Then there is $u\in H^{1}(\R^{3})\backslash\{0\}$
with $u\geq 0$ such that, passing a subsequence, we have $u_{n}\rig
u$ in $H^{1}(\R^{3})$. Moreover, $u$ is a positive ground-state
solution for the problem $(P_{\infty})$.

\end{Th}

\noindent {\bf Proof.} Firstly we prove that $(u_{n})$ is bounded in
$H^{1}(\R^{3})$. Since $(u_{n})\subset H^{1}(\R^{3})$ is  a
$(PS)_{c_{\infty}}$ sequence for $I_{\infty}$, there exists  $C>0$
such that
$$
C+\|u_{n}\|\geq
I_{\infty}(u_{n})-\frac{1}{\theta}I_{\infty}'(u_{n})u_{n}, \ \forall
n\in \N.
$$
From $(f_{3})$, we get
$$
C+\|u_{n}\|\geq\frac{1}{2}\cha{M}(\|u_{n}\|^{2})-\frac{1}{\theta}M(\|u_{n}\|^{2})\|u_{n}\|^{2},
\ \forall n\in \N.
$$
Notice that by $(M_{3})$ we have
\begin{eqnarray}\label{permitelimitacaoPS}
\widehat{M}(t)\geq \frac{1}{2}M(t)t,
\end{eqnarray}
for all $t\geq 0$. This inequality allows us to conclude that,
\begin{eqnarray}\label{permitelimitacaoPS}
\frac{1}{t}\left[
\frac{1}{2}\cha{M}(t)-\frac{1}{\theta}M(t)t\right]\geq\left(\frac{\theta
-4}{4\theta}\right)m_{0},
\end{eqnarray}
for all $t>0$. Now, suppose by contradiction  that, up to a
subsequence, $\|u_{n}\|\rig \infty$. Thus,
$$
\frac{C}{\|u_{n}\|^{2}}+\frac{1}{\|u_{n}\|}\geq\frac{1}{\|u_{n}\|^{2}}
\left[\frac{1}{2}\cha{M}(\|u_{n}\|^{2})-\frac{1}{\theta}M(\|u_{n}\|^{2})\|u_{n}\|^{2}\right],
\ \forall n\in \N.
$$
Passing to limit we get
$$
0\geq \left(\frac{\theta-4}{4\theta}\right)m_{0}>0,
$$
which is an absurd. Hence, there exist $u \in H^{1}(\mathbb{R}^{3})$
and $t_{0}>0$ such that
\begin{equation}
u_{n}\righ u \ \mbox{in} \ H^{1}(\R^{3}),\label{weak convergence 1}
\end{equation}
and
\begin{equation}
\|u_{n}\|\rig t_{0}.\label{norm convergence 1}
\end{equation}

Since $M$ is continuous function, we get
\begin{equation}\label{JRSJ}
M(\|u_{n}\|^{2})\rig M(t_{0}^{2}).
\end{equation}

From \cite[Theorem 2.5]{Alves and Figueiredo} we conclude that
\begin{equation}\label{JRSJ1}
M(t_{0}^{2})=M(\|u\|^{2}).
\end{equation}

Using $(M_{2})$ we obtain $\|u\|=t_{0}$ and the lemma is proved. $\fim$

\vspace*{.3cm}

An important property that we can derive from (\ref{weak convergence 1}), (\ref{norm convergence 1}) and
(\ref{JRSJ1}) is that, up to a subsequence,
\begin{equation}\label{PS}
u_{n}\rightarrow u \ \ \mbox{in} \ \ H^{1}(\mathbb{R}^{3})
\end{equation}

\subsection{A compactness result on the Nehari manifold associated to limit problem}

We denote by $H^{1, +}(\R^{3})$ the open subset of $H^{1}(\R^{3})$
given by
$$
H^{1,+}(\R^{3})=\{u\in H^{1}(\R^{3}): u^{+}=\max\{0,u\}\neq 0\},
$$
and $S_{\infty}^{+}=S_{\infty}\cap H^{1,+}(\R^{3})$, where
$S_{\infty}$ is unit sphere of $H^{1}(\R^{3})$.

Note that $S_{\infty}^{+}$ is a non-complete
$C^{1,1}$-manifold of codimension $1$, modeled on $H^{1}(\R^{3})$ and contained  in the open
$H^{1, +}(\R^{3})$. Thus, $H^{1}(\R^{3})=T_{u}S_{\infty}^{+}\oplus
\R \ u$ for each $u\in S_{\infty}^{+}$, where
$T_{u}S_{\infty}^{+}=\{v\in H^{1}(\R^{3}):(u, v)=0\}$.

As can be seen in \cite[Theorem 3.1]{Alv1}, to prove the main
result, it is very important to obtain a result of compactness on of
the Nehari manifold associated to limit problem. Since that $f$ and
$M$ are only continuous, we cannot claim that $\mathcal{N}_{\infty}$
is a continuous manifold. Here was necessary a new argument. In the
Lemmas \ref{lema3.2.2} and \ref{Proposition 3.1.2} we adapt
arguments from the excellent book \cite[Chapter 3]{Szulkin}, see also \cite{Szulkin 1}. In the
Lemma \ref{PVE} we prove a result of the type Ekeland's Principle on
set $S^{+}_{\infty}$. The behavior of the functional $I_{\infty}$
near the boundary of $S^{+}_{\infty}$ (see Lemma \ref{lema3.2.2} $(A_{4})$) overcome the fact that this set is a
non-complete $C^{1,1}$-manifold of $H^{1}(\R^{3})$. In the
Lemma \ref{compacityLemma} we prove the main result of the section.

\begin{Le}\label{lema3.2.2}
Suppose that the function $M$ satisfies $(M_{1})-(M_{3})$  and the
function $f$ satisfies $(f_{1})-(f_{4})$. Then:

\begin{description}
\item[($A_{1}$)] For each $u\in H^{1,+}(\mathbb{R}^{3})$, let
$h:\R_{+}\rig\R$ be defined by $h_{u}(t)=I_{\infty}(tu)$. Then,
there is a unique $t_{u}>0$ such that $h_{u}'(t)>0$ in $(0,t_{u})$
and $h_{u}'(t)<0$ in $(t_{u}, \infty)$.

\item[($A_{2}$)] there is $\tau>0$ independent on $u$ such that $t_{u}\geq \tau$ for all $u\in S^{+}_{\infty}$. Moreover,
for each compact set $\mathcal{W}\subset S^{+}_{\infty}$ there is
$C_{\mathcal{W}}>0$ such that $t_{u}\leq C_{\mathcal{W}}$, for all
$u\in \mathcal{W}$.\label{Lemma 4.1.2}

\item[($A_{3}$)] The map
$\widehat{m}_{\infty}:H^{1,+}(\mathbb{R}^{3})\rightarrow
\mathcal{N}_{\infty}$ given by $\widehat{m}_{\infty}(u)=t_{u}u$ is
continuous and
$m_{\infty}:=\widehat{m}_{\infty_{\bigl|S^{+}_{\infty}}}$ is a
homeomorphism between $S^{+}_{\infty}$ and $\mathcal{N}_{\infty}$.
Moreover, $m_{\infty}^{-1}(u)= \frac{u}{\|u\|}$.

\item[($A_{4}$)] If there is a sequence $(u_{n})\subset
S^{+}_{\infty}$ such that \mbox{dist}$(u_{n},\partial
H^{1,+}(\mathbb{R}^{3}))\rightarrow 0$, then
$\|m_{\infty}(u_{n})\|\rightarrow \infty$ and
$I_{\infty}(m_{\infty}(u_{n}))\rightarrow \infty$.
\end{description}
\end{Le}
\noindent {\bf Proof.} Since that $(M_{1})$ and (\ref{crescimentoM})
occurs, the items ($A_{1}$), ($A_{2}$) and ($A_{3}$) follows of the a
simple adaptation of the \cite[Proposition 8]{Szulkin}. The item
($A_{4}$) follows by \cite[Lemma 26]{Szulkin}. ${\fim}$

We use $(A_{3})$ and $(A_{4})$ to overcome the lack of
differentiability of $\mathcal{N}_{\infty}$. (see Lemmas \ref{PVE}
and \ref{compacityLemma}).

Now we define
$$
\cha{\Psi}_{\infty}:H^{1, +}(\R^{3})\rig\R \ \mbox{e} \
\Psi_{\infty}:S_{\infty}^{+}\rig \R,
$$
by $\cha{\Psi}_{\infty}(u)=I_{\infty}(\cha{m}_{\infty}(u))$ \ and
$\Psi_{\infty}:=(\cha{\Psi}_{\infty})_{|_{S_{\infty}^{+}}}$.

Using the same type of arguments explored in (\cite[Lemma
26]{Szulkin}) we can prove the next lemma. Thus, we omit the proof.

\begin{Le}\label{Proposition 3.1.2}
Suppose that the function $M$ satisfies $(M_{1})-(M_{3})$  and the
function $f$ satisfies $(f_{1})-(f_{4})$. Then:
\begin{description}
\item[($a$)] $\cha{\Psi}_{\infty}\in C^{1}(H^{1, +}(\R^{3}), \R)$ e
$$
\cha{\Psi}_{\infty}'(u)v=\frac{\|\cha{m}_{\infty}(u)\|}{\|u\|}I_{\infty}'(\cha{m}_{\infty}(u))v,
\ \forall u\in H^{1, +}(\R^{3}) \ \mbox{e} \ \forall v\in
H^{1}(\R^{3}).
$$

\item[($b$)] $\Psi_{\infty}\in C^{1}(S_{\infty}^{+}, \R)$ e
$$
\Psi_{\infty}'(u)v=\|m(u)\|I_{\infty}'(m_{\infty}(u))v, \ \forall
v\in T_{u}S_{\infty}^{+}.
$$

\item[($c$)] If  $(u_{n})$ is a  $(PS)_{c}$ sequence for $\Psi_{\infty}$, then $(m_{\infty}(u_{n}))$ is a  $(PS)_{c}$ sequence for the
functional $I_{\infty}$. If $(u_{n})\subset \mathcal{N}_{\infty}$ is
a  bounded $(PS)_{c}$ sequence for $I_{\infty}$, then
$(m^{-1}_{\infty}(u_{n}))$ is a  $(PS)_{c}$ sequence for
$\Psi_{\infty}$.

\item[($d$)] $u$ is critical point of $\Psi_{\infty}$ if, and only if, $m_{\infty}(u)$ is a nontrivial critical point of $I_{\infty}$.
Moreover, the critical values are the same and
$$
\inf_{S_{\infty}^{+}}\Psi_{\infty}=\inf_{\mathcal{N}_{\infty}}I_{\infty}.
$$
\end{description}
\end{Le}

By using $(M_{1})-(M_{3})$ we have, as in \cite[Remark 11, Remark 34]{Szulkin},
the following variational characterization of the infimum of
$I_{\infty}$ over $\mathcal{N}_{\infty}$:

\begin{eqnarray}\label{minimax}
c_{\infty}=\inf_{u\in \mathcal{N}_{\infty}}I_{\infty}(u)=\inf_{u\in
H^{1, +}(\R^{3})}\max_{t>0}I_{\infty}(tu)=\inf_{u\in
S_{\infty}^{+}}\max_{t>0}I_{\infty}(tu).
\end{eqnarray}

\begin{Le}\label{PVE}
Let $(v_{n}) \subset S_{\infty}^{+}$ be a sequence such that
$\Psi_{\infty}(v_{n})\rig c_{\infty}$. Then, there is a sequence
$(\widehat{v}_{n})\subset S_{\infty}^{+}$ such that
$(\widehat{v}_{n})$ is a  $(PS)_{c_{\infty}}$sequence for
$\Psi_{\infty}$ in $S_{\infty}^{+}$ and
$\|\widehat{v}_{n}-v_{n}\|=o_{n}(1)$.
\end{Le}

\noindent {\bf Proof.} Let $(V, d)$ be a
complete metric space, where $V=\overline{H^{1,+}(\R^{3})}$ and $d(u, v)=\|u-v\|$, 
and a map  $\zeta_{\infty}:V\rig
\R\cup\{\infty\}$ given by $\zeta(u)=\cha{\Psi}_{\infty}(u)$ if
$u\in H^{1,+}(\R^{3})$ and $\zeta(u)=\infty$ if $u\in \partial
H^{1,+}(\R^{3})$. From Lemma  $\ref{lema3.2.2}$ $(A_{4})$, we have
that $\zeta$ is continuous and by $(M_{3})$ and $(f_{3})$ we
conclude that $\zeta$ is bounded below. Using Eleland's Variational
Principle \cite[Theorem 1.1]{Ekeland}, it follows that for each given $\var, \lambda>0$ and each $u\in V$, with
$c_{\infty}<\zeta(u)< c_{\infty}+\var$, there is $v\in V$ such that
\begin{equation}\label{114}
\zeta(v)\leq \zeta(u), \|u-v\|\leq \lambda \ \mbox{e} \
\zeta(w)>\zeta(v)-(\frac{\var}{\lambda})\|v-w\|, \ \forall \ w\neq
v.
\end{equation}
Once that $(v_{n})\subset S_{\infty}^{+}$ and
$\zeta(v_{n})=\Psi_{\infty}(v_{n})\rig c_{\infty}$, from Lemma
$\ref{lema3.2.2}$ $(A_{4})$, there exists $R>0$ which independent on
$n\in\N$ such that
$$
dist(v_{n}, \partial H^{1, +}(\R^{3}))>R, \ \forall \ n\in\N.
$$

Thus, for each  $z\in \overline{B_{1}(0)}\subset H^{1}(\R^{3})$ and
for each $n\in\N$, we get
\begin{equation}\label{118}
v_{n}+tz\in B_{ \frac{R}{2}}(v_{n})\subset H^{1, +}(\R^{3}), \
\forall \ t\in(0, \frac{R}{2}).
\end{equation}

We can consider without loss of generality  that
$\Psi_{\infty}(v_{n})<c_{\infty}+\dis\frac{1}{n}$. Hence, choosing
$u=v_{n}$ and $\var=\lambda=\frac{1}{n}$ in $(\ref{114})$, we obtain
$\widehat{w}_{n}\in V$ such that
\begin{equation}\label{117}
\zeta_{\infty}(\widehat{w}_{n})\leq\Psi_{\infty}(v_{n}),
\end{equation}
\begin{equation}\label{115}
\|\widehat{w}_{n}-v_{n}\|\leq \frac{1}{n}
\end{equation}
and
\begin{equation}\label{116}
\zeta_{\infty}(\widehat{w}_{n}+tz)>\Psi_{\infty}(\widehat{w}_{n})-(\frac{1}{n^{2}})\|tz\|,
\ \forall \ t\in(0, \frac{R}{2}),
\end{equation}
where in $(\ref{114})$ was chosen $w=\widehat{w}_{n}+tz$. From
$(\ref{117})$, $(\ref{115})$ and  $(\ref{118})$ we derive
$(\widehat{w}_{n})\subset H^{1, +}(\R^{3})$ and
$$
\|\widehat{w}_{n}+tz-v_{n}\|< \frac{1}{n}+\frac{R}{2}, \ \forall \
n\in\N \ \mbox{e} \ \forall \ t\in(0, \frac{R}{2}).
$$

Thus, for $n$ large and $t\in(0, \frac{R}{2})$, we have
$\widehat{w}_{n}+tz\in B_{ \frac{R}{2}}(v_{n})\subset H^{1,
+}(\R^{3})$ and consequently
$\zeta_{\infty}(\widehat{w}_{n})=\cha{\Psi}_{\infty}(\widehat{w}_{n})$
and
$\zeta_{\infty}(\widehat{w}_{n}+tz)=\cha{\Psi}_{\infty}(\widehat{v}_{n}+tz)$.

By $(\ref{116})$ we obtain
\begin{equation}\label{121}
\frac{\cha{\Psi}_{\infty}(\widehat{w}_{n}+tz)-\cha{\Psi}_{\infty}(\widehat{w}_{n})}{t}>-(\frac{1}{n^{2}})\|z\|,
\end{equation}
for all $t\in (0, \frac{R}{2}), \ n\geq n_{0} \ \mbox{and} \ z\in
\overline{B_{1}(0)}$. Using the definition of $\cha{\Psi}_{\infty}$
follow that, for each $u\in H^{1, +}(\R^{3})$, we get
$\cha{\Psi}_{\infty}(t u)=\cha{\Psi}_{\infty}(u)$, for all $t>0$.
Now defining
$\cha{v}_{n}=\frac{\widehat{w}_{n}}{\|\widehat{w}_{n}\|}$, we
conclude that $\{\cha{v}_{n}\}\subset S_{\infty}^{+}$. Moreover,
from $(\ref{117})$ we derive
\begin{equation}\label{119}
c_{\infty}\leq
\Psi_{\infty}(\widehat{v}_{n})\leq\Psi_{\infty}(v_{n}).
\end{equation}
By $(\ref{115})$, we obtain
\begin{equation}\label{122}
1-\frac{1}{n}<\|\widehat{w}_{n}\|<1+\frac{1}{n}, \ \forall \ n\in\N.
\end{equation}
From $(\ref{115})$ and a straightforward computation, we have
\begin{equation}\label{120}
\| \widehat{v}_{n}-v_{n}\|\leq \frac{2}{n-1}
\end{equation}
Finally, from $(\ref{121})$ we conclude
$$
\frac{\cha{\Psi}_{\infty}(\frac{\widehat{w}_{n}+tz}{\|\widehat{w}_{n}\|})-\Psi_{\infty}(\widehat{v}_{n})}{t}>-(\frac{1}{n^{2}})\|z\|,
$$
for all $t\in (0, \frac{R}{2}), \ n\geq n_{0} \ \mbox{ande} \ z\in
\overline{B_{1}(0)}$, that implies
$$
\frac{\cha{\Psi}_{\infty}(\widehat{v}_{n}+\frac{t}{\|\widehat{w}_{n}\|}z)-\Psi_{\infty}(\widehat{v}_{n})}{\frac{t}{\|\widehat{w}_{n}\|}}>-(\frac{1}{n^{2}})\|z\|\|\widehat{w}_{n}\|,
$$
for all  $t\in (0, \frac{R}{2}), \ n\geq n_{0} \ \mbox{e} \ z\in
\overline{B_{1}(0)}$. Passing to the limit in $t\rig 0$ and using
$(\ref{122})$ we have
$$
\cha{\Psi}'_{\infty}(\widehat{v}_{n})z\geq
-(\frac{1}{n^{2}})(1+\frac{1}{n})\|z\|
$$
for all $\ n\geq n_{0} \ \mbox{and} \ z\in \overline{B_{1}(0)}$.
Considering now $\ z\in \overline{B_{1}(0)}\cap
T_{\widehat{v}_{n}}S_{\infty}$, we get
\begin{equation}\label{123}
\|\Psi'_{\infty}(\widehat{v}_{n})\|_{\ast}\leq
(\frac{1}{n^{2}})(1+\frac{1}{n}),
\end{equation}
for all $\ n\geq n_{0}$. Passing to the limit in $n\rig\infty$
in $(\ref{119})$, $(\ref{120})$ and $(\ref{123})$ we prove the
Lemma. ${\fim}$.

The following compactness property will be crucial in our arguments
and is necessary because the Nehari manifold is not a regular
manifold. If $M$ and $f$ are $C^{1}$ functions, we can to argue as in \cite[Proposition 3.1]{Alv1}.

\begin{Le}\label{compacityLemma}
Let $(u_{n}) \subset \mathcal{N}_{\infty}$ be a sequence such that
$I_{\infty}(u_{n})\rig c_{\infty}$. Then,
$$
u_{n}(x)=w_{n}(x-y_{n})+\Psi(x-y_{n}),
$$
where $\{w_{n}\}\subset H^{1}(\R^{3})$ with
$$
w_{n}\rig 0, \ \mbox{em} \ H^{1}(\R^{3}),
$$
$\{ y_{n}\}\subset \R^{3}$ is such that $|y_{n}|\rig \infty$ and
$\Psi \in H^{1}(\R^{3})$ is a positive continuous function  that
satisfies
$$
I_{\infty}(\Psi)=c_{\infty} \ \mbox{e} \ I_{\infty}'(\Psi)\Psi=0.
$$
\end{Le}

\noindent {\bf Proof.} Follow from Lemmas  $\ref{lema3.2.2}$ ($(A_{3})$),
 $\ref{Proposition 3.1.2}$ and
$\ref{minimax}$, that
\begin{equation}\label{equality 1}
v_{n}=m^{-1}_{\infty}(u_{n})=\frac{u_{n}}{\|u_{n}\|}\in
S_{\infty}^{+}, \ \forall n\in \N
\end{equation}
and
$$
\Psi_{\infty}(v_{n})=I_{\infty}(u_{n})\rig
c_{\infty}=\inf_{S_{\infty}^{+}}\Psi_{\infty}.
$$

By Lemma $\ref{PVE}$, there is a sequence
$\{\widehat{v}_{n}\}\subset S_{\infty}^{+}$ such that
$\{\widehat{v}_{n}\}$ is a  $(PS)_{c_{\infty}}$ sequence for
$\Psi_{\infty}$ on $S_{\infty}^{+}$ and
\begin{equation}\label{convergencia}
\|\widehat{v}_{n}-v_{n}\|_{\infty}=o_{n}(1).
\end{equation}

Thus,
$$
\cha{\Psi}_{\infty}'(\widehat{v}_{n})v=2\lambda_{n}(\widehat{v}_{n},
v)+o_{n}(1), \ \forall v\in H^{1}(\R^{3}).
$$
From Lemma $\ref{Proposition 3.1.2}$ ($(a)$), we have
$$
\|m_{\infty}(\widehat{v}_{n})\|I_{\infty}'(m_{\infty}(\widehat{v}_{n}))v=2\lambda_{n}(\widehat{v}_{n},
v)+o_{n}(1), \ \forall v\in H^{1}(\R^{3}),
$$
that implies
\begin{equation}
I_{\infty}'(m_{\infty}(\widehat{v}_{n}))v=2\lambda_{n}\left(\widehat{v}_{n},
\frac{v}{\|m_{\infty}(\widehat{v}_{n})\|}\right)+o_{n}(1), \ \forall
v\in H^{1}(\R^{3}).\label{equation 9}
\end{equation}

Hence, for $v=m_{\infty}(\widehat{v}_{n})$, we get
$$
0=I_{\infty}'(m_{\infty}(\widehat{v}_{n}))m_{\infty}(\widehat{v}_{n})=2\lambda_{n}+o_{n}(1).
$$

Thus, we conclude that $\lambda_{n}\rig 0$ as $n\rig \infty$ and
from $(\ref{equation 9})$, we obtain
$$
\|I_{\infty}'(m_{\infty}(\widehat{v}_{n}))\|\leq
C|\lambda_{n}|+o_{n}(1)=o_{n}(1).
$$

Thus, $\{m_{\infty}(\widehat{v}_{n})\}$ is a $(PS)_{c_{\infty}}$
sequence for $I_{\infty}$ in $H^{1}(\R^{3})$. From  $(\ref{PS})$,
we obtain $\{\cha{w}_{n}\}\subset H^{1}(\R^{3})$ and $\{ y_{n}\}\subset \R^{3}$ such that 
$$
\cha{w}_{n}\rig 0, \ \mbox{in} \ H^{1}(\R^{3}),
$$
$|y_{n}|\rig \infty$ and
$$
m_{\infty}(\widehat{v}_{n})(x)=\cha{w}_{n}(x-y_{n})+\cha{\Psi}(x-y_{n}),
$$
where $\cha{\Psi} \in H^{1}(\R^{3})$ is a positive function satisfying
\begin{equation}\label{variedade}
I_{\infty}(\cha{\Psi})=c_{\infty} \ \mbox{and} \ I_{\infty}'(\cha{\Psi})\cha{\Psi}=0.
\end{equation}
So $\|\cha{\Psi}\|\geq \tau>0$, and defining
$$
\widehat{v}_{n}(x)=\frac{m_{\infty}(\widehat{v}_{n})(x)}{\|m_{\infty}(\widehat{v}_{n})\|}=
\frac{\cha{w}_{n}(x-y_{n})}{\|m_{\infty}(\widehat{v}_{n})\|}+ \frac{\cha{\Psi}(x-y_{n})}{\|m_{\infty}(\widehat{v}_{n})\|},
$$
we conclude that
$$
\widehat{v}_{n}(x)=w_{n}(x-y_{n})+\dis\frac{\cha{\Psi}(x-y_{n})}{\|m_{\infty}(\widehat{v}_{n})\|},
$$
where $w_{n}(x-y_{n}):= \dis\frac{\cha{w}_{n}(x-y_{n})}{\|m_{\infty}(\widehat{v}_{n})\|}$ is such that $w_{n}\rig 0$ in $H^{1}(\R^{3})$.
Follow from (\ref{convergencia}) that
$$
v_{n}(x)=w_{n}(x-y_{n})+\dis\frac{\cha{\Psi}(x-y_{n})}{\|m_{\infty}(\widehat{v}_{n})\|}.
$$
From $(\ref{equality 1})$ and from the continuity of $m_{\infty}$, we have 
$$
u_{n}(x)=w_{n}(x-y_{n})+m_{\infty}\left(\dis\frac{\cha{\Psi}(x-y_{n})}{\|m_{\infty}(\widehat{v}_{n})\|}\right).
$$
Defining, $\Psi:= m_{\infty}\left(\dis\frac{\cha{\Psi}}{\|m_{\infty}(\widehat{v}_{n})\|}\right)=m_{\infty}(\cha{\Psi})$, 
it follows that
$$
u_{n}(x)=w_{n}(x-y_{n})+\Psi(x-y_{n}),
$$
with $w_{n}\rig 0$ in  $H^{1}(\R^{3})$. From $(\ref{variedade})$, we derive
$$
I_{\infty}(\Psi)=c_{\infty} \ \mbox{e} \ I_{\infty}'(\Psi)\Psi=0.
$$

 $\fim$

\section{Variational framework and behavior of minimax levels}

From now on we will assume, without loss of generality, that $0\in
\Omega$. Let us fix real numbers $R>r>0$ such that $B_r(0) \subset
\Omega \subset B_R(0)$ and the sets
$$
\Omega^{+}:=\{x \in \mathbb{R}^N :
\mbox{dist}(x,\overline{\Omega})\leq r\},~~~ \Omega^{-}:=\{x\in
\Omega : \mbox{dist}(x,\partial\Omega)\geq r\}
$$
are homotopically equivalent to $\Omega$.

For each $\lambda>0$, we shall denote by
$H^{1}_{0}(\Omega_{\lambda})$ the Hilbert space obtained by the
closure of $C^{\infty}_{0}(\Omega_{\lambda})$ under the scalar
product
$$
\langle u,v \rangle_{\lambda} := \int_{\Omega_{\lambda}} \nabla u
\nabla v \ dx + \int_{\Omega_{\lambda}}uv \ dx.
$$

The norm induced by this inner product is given by
$$
\|u\|_{\lambda}:= \left(\int_{\Omega_{\lambda}} |\nabla u|^2 \ dx+
\int_{\Omega_{\lambda}}|u|^2 \ dx\right)^{1/2}.
$$

In view of $(f_1)-(f_2)$, we have that the functional
$I_{\lambda}:H^{1}_{0}(\Omega_{\lambda}) \to \mathbb{R}$ given by
\begin{equation} \label{def1}
J_{\lambda}(u) :=
\displaystyle\frac{1}{2}\widehat{M}(\|u\|^{2}_{\lambda})
 -\displaystyle\int_{\Omega_{\lambda}} F(u) \ dx
\end{equation}
is well defined. Moreover, $J_{\lambda}  \in
C^1(H^{1}_{0}(\Omega_{\lambda}))$ with the following derivative
$$
J'_{\lambda}(u)v=M(\|u\|^{2}_{\lambda})\biggl[\displaystyle\int_{\Omega_{\lambda}}\nabla
u\nabla v \ dx + \displaystyle\int_{\Omega_{\lambda}} u v \
dx\biggl] -\displaystyle\int_{\Omega_{\lambda}} f(u)v \ dx.
$$
Thus the weak solutions of $(P_{\lambda})$ are  precisely the
critical points of $J_{\lambda}$.

As in the previous section, we denote by $H^{1,
+}(\Omega_{\lambda})$ the open subset of
$H^{1}_{0}(\Omega_{\lambda})$ given by
$$
H^{1,+}(\Omega_{\lambda})=\{u\in H^{1}_{0}(\Omega_{\lambda}):
u^{+}=\max\{0,u\}\neq 0\},
$$
and $S_{\lambda}^{+}=S_{\lambda}\cap H^{1,+}(\Omega_{\lambda})$,
where $S_{\lambda}$ is unit sphere of $H^{1}_{0}(\Omega_{\lambda})$.
Recalling that
$H^{1}_{0}(\Omega_{\lambda})=T_{u}S_{\lambda}^{+}\oplus \R \ u$ for
each $u\in S_{\lambda}^{+}$, where $T_{u}S_{\lambda}^{+}=\{v\in
H^{1}_{0}(\Omega_{\lambda}): \langle u, v\rangle_{\lambda}=0\}$.

In view of the subcritical growth of $f$, $(f_3)$, $(M_{1})$ and
(\ref{crescimentoM}), it is standard to check that $J_{\lambda}$
satisfies the Palais-Smale condition. Moreover, these hypotheses
imply that $J_{\lambda}$ has the mountain pass geometry. Hence, for
each $\lambda>0$, there exists $u_{\lambda} \in
H^{1}_{0}(\Omega_{\lambda})$ such that
$J_{\lambda}(u_{\lambda})=b_{\lambda}$ and
$J_{\lambda}'(u_{\lambda})=0$, where $b_{\lambda}$ denotes the
mountain pass level of the functional $J_{\lambda}$. 

\begin{Remark}\label{remark}
We point out that analogous results to the Lemmas \ref{lema3.2.2} and 
\ref{Proposition 3.1.2} are still true for the functional $J_{\lambda}$.
Moreover, in the present case, we denote the functions of those Lemmas by 
$\cha{m}_{\lambda}, m_{\lambda}, \cha{\Psi}_{\lambda}$ and $\Psi_{\lambda}$.
\end{Remark}

By using
$(M_{1})-(M_{3})$ and as \cite[Remark 11, Remark 34]{Szulkin},  we
can prove that $b_{\lambda}$ can also be characterized as
\begin{equation} \label{def2}
b_{\lambda}=\inf_{u\in
\mathcal{M}_{\lambda}}J_{\lambda}(u)=\inf_{u\in H^{1,
+}(\Omega_{\lambda})}\max_{t>0}J_{\lambda}(tu)=\inf_{u\in
S_{\lambda}^{+}}\max_{t>0}J_{\lambda}(tu),
\end{equation}
where $\mathcal{M}_{\lambda}$ is the Nehari manifold associated to
$J_{\lambda}$, namely
\begin{equation} \label{def3}
\mathcal{M}_{\lambda} := \{u \in H^{1}_{0}(\Omega_{\lambda})
\backslash\{0\}: J_{\lambda}'(u)u=0\}.
\end{equation}
From $(f_1)$ and $(f_2)$, there exists $r=r(\lambda)>0$ such that
\begin{equation}
\label{prop-nehari} \|u\|_{\lambda} \geq r >0,
\end{equation}
for all $u \in \mathcal{M}_{\lambda}$.

We recall that $B_{\lambda r}(0) \subset \Omega_{\lambda}$ and
define the triple
$(J_{\lambda,r},\,b_{\lambda,r},\,\mathcal{M}_{\lambda,r})$  in a
similar way, just replacing $\Omega_{\lambda}$ by $B_{\lambda
r}(0)$.

For each $x\in \mathbb{R}^3$, let us denote by $A_{\lambda,x}$ the
following set
$$
A_{\lambda,x}:=B_{\lambda R}(x)\setminus \overline{B_{\lambda r}(x)}
$$
and define the functional
$\widehat{J}_{\lambda,x}:H^{1}_{0}(A_{\lambda ,x}) \to \mathbb{R}$
by
\begin{equation}
\widehat{J}_{\lambda,x}(v):=
\frac{1}{2}\widehat{M}\biggl(\int_{A_{\lambda,x}} |\nabla v |^2 \ dx
+ \frac{1}{2}\int_{A_{\lambda ,x}} |v|^2\ dx \biggl)  -
\int_{A_{\lambda ,x}} F(v) \ dx.
\end{equation}
and the set
$$
\widehat{\mathcal{A}}_{\lambda, x}:=\{v \in H^{1}_{0}(A_{\lambda
,x}, \mathbb{R})\backslash \{0\}: \widehat{J}'_{\lambda,x}(v)v=0 \}.
$$

For $v \in H^{1}(\mathbb{R}^{N})$ with compact support, we consider
the barycenter map
$$
\beta(v):=\frac{\displaystyle\int_{\mathbb{R}^3}x|\nabla v|^{2}\
dx}{\displaystyle\int_{\mathbb{R}^3}|\nabla v|^{2}\ dx}
$$
and we introduce the following quantity
$$
a_{\lambda,x}:=\inf\left\{\widehat{J}_{\lambda,x}(v):
v\in\widehat{{\cal{A}}}_{\lambda,x}\;\;\mbox{and}\;\; \beta(v)=x
\right\}.
$$

We present below an important property of the asymptotic behavior of
the numbers $a_{\lambda,0}$.

\begin{Le} \label{lema-assintotico2}
The following holds
$$
c_{\infty} < \liminf_{\lambda \to\infty} a_{\lambda,0}.
$$
\end{Le}

\noindent {\bf Proof.} Since $c_{\infty} \leq a_{\lambda,0}$ for any
$\lambda>0$, we have that $c_{\infty} \leq \liminf_{\lambda
\to\infty} a_{\lambda,0}$. Suppose, by contradiction, that for some
sequence $\lambda_n \nearrow \infty$ we have that $a_{\lambda_n,0}
\to c_{\infty}$. Then, we can obtain $v_n \in \widehat{\mathcal{A}}_{\lambda_n,0}
\subset \mathcal{N}_{\infty}$ satisfying
$\widehat{J}_{\lambda_n,0}(v_n) = I_{\infty}(v_n) \to c_{\infty}$
and $\beta(v_n)=0$, where we are understanding that the function
$v_n$ is extended to the whole space by setting $v_n(x):=0$ for a.e.
$x \in \mathbb{R}^N \setminus A_{\lambda_n,0}$.

Thus, it follows from
Lemma \ref{compacityLemma} that
\begin{equation} \label{sono3}
v_n(x) = w_n(x-y_n) + \widetilde{v}(x-y_n)
\end{equation}
with $(w_n) \subset  H^1(\mathbb{R}^3)$ satisfying $w_n \to 0$
strongly in $H^1(\mathbb{R}^3)$, $(y_n) \subset \mathbb{R}^3$ being
such that $|y_n|\to\infty$, and $\widetilde{v} \in
H^1(\mathbb{R}^3)$ verifying
\begin{equation} \label{sol_infinito}
J_{\infty}(\widetilde{v}) =
c_{\infty},~~~J_{\infty}'(\widetilde{v})=0.
\end{equation}

The rest of the proof follows as in \cite[Proposition 4.1]{Alv1}.
${\fim}$

In the next result, we present the asymptotic behavior of the
minimax $b_{\lambda}$ as $\lambda \to \infty$.

\begin{Le} \label{lema-assintotico1}
We have that
$$
\lim_{\lambda \to \infty} b_{\lambda}=c_{\infty}.
$$
\end{Le}
\noindent {\bf Proof.} Let $\varphi\in C_{0}^{\infty}(\R^{3})$ be a
function such that $\varphi=1$ in $B_{1}(0)$, $\varphi=0$ in
$\R^{3}\backslash B_{2}(0)$ and $0\leq \varphi\leq 1$. For each
$R>0$, we define
$$
\varphi_{R}(x)=\varphi\left(\frac{x}{R}\right) \ \mbox{and} \
w_{R}(x)=\varphi_{R}(x)w(x),
$$
where $w$ is a ground-state solution of $(P_{\infty})$. Arguing as
\cite[Proposition 4.2]{Alv1} we conclude that
\begin{equation}\label{40}
\limsup_{\lambda\to \infty}c_{\lambda}\leq I_{\infty}(t_{R}w_{R}).
\end{equation}

Now we show that
$$
\dis\lim_{R\to\infty}t_{R}=1.
$$

Indeed, since  $\|w_{R}\|_{\lambda}^{2}=\|w_{R}\|^{2}$, we have
\begin{equation}\label{41}
\frac{M(t_{R}^{2}\|w_{R}\|^{2})}{t_{R}^{2}\|w_{R}\|^{2}}=
\frac{1}{\|w_{R}\|^{4}}\dis\int_{\R^{3}}\left[
\frac{f(t_{R}w_{R})}{(t_{R}w_{R})^{3}}\right]w_{R}^{4}.
\end{equation}

From $(f_{3})$ e $(f_{4})$ and $R>1$, we get
\begin{equation}\label{42}
\frac{M(t_{R}^{2}\|w_{R}\|^{2})}{t_{R}^{2}\|w_{R}\|^{2}}\geq
\frac{1}{\|w_{R}\|^{4}}\dis\int_{B_{1}(0)}\left[\frac{f(t_{R}w)}{(t_{R}w)^{3}}\right]w^{4}\geq
\frac{1}{\|w_{R}\|^{4}}\dis\int_{B_{1}(0)}\left[\frac{f(t_{R}a)}{(t_{R}a)^{3}}\right]a^{4},
\end{equation}
where $a=\dis\min_{|x|\leq 1}w(x)$.

Suppose that there is $(R_{n})$ a sequence such that $t_{R_{n}}\rig
\infty$ as $R_{n}\rig\infty$. Thus, we have
$$
\frac{M(t_{R_{n}}^{2}\|w_{R_{n}}\|^{2})}{t_{R_{n}}^{2}\|w_{R_{n}}\|^{2}}\geq
\frac{1}{\|w_{R_{n}}\|^{4}}\dis\int_{B_{1}(0)}\left[\frac{f(t_{R_{n}}a)}{(t_{R_{n}}a)^{3}}\right]a^{4}.
$$
It follows from $(f_{3})$ and Fatou's Lemma that
\begin{equation*}
\frac{M(t_{R_{n}}^{2}\|w_{R_{n}}\|^{2})}{t_{R_{n}}^{2}\|w_{R_{n}}\|^{2}}\rig\infty,
\end{equation*}
which is a contradiction with $(M_{3})$.

Suppose now $t_{R_{n}}\rig 0$ as $R_{n}\rig\infty$. Using the growth
of $f$, given by $(f_{1})-(f_{2})$, we have
\begin{equation*}
\frac{M(t_{R_{n}}^{2}\|w_{R_{n}}\|^{2})}{t_{R_{n}}^{2}\|w_{R_{n}}\|^{2}}\rig
0,
\end{equation*}
which is a contradiction with $(M_{1})$. Thus, there exists $t_{0}>0$
such that, up to a subsequence, $t_{R_{n}}\rig t_{0}$ and from
$(\ref{41})$, we get
\begin{equation*}
\frac{M(t_{0}^{2}\|w\|^{2})}{t_{0}^{2}\|w\|^{2}}=
\frac{1}{\|w\|^{4}}\dis\int_{\R^{3}}\left[
\frac{f(t_{0}w)}{(t_{0}w)^{3}}\right]w^{4}.
\end{equation*}

Since that $w$ is a solution of $(P_{\infty})$ we conclude that
$t_{0}=1$ and $I_{\infty}(t_{R}w_{R})\rig I_{\infty}(w)=c_{\infty}$
as $R\rig \infty$. Thus, by $(\ref{40})$ we obtain
\begin{equation*}
\limsup_{\lambda\to \infty}c_{\lambda}\leq c_{\infty}.
\end{equation*}
The reverse inequality  follows from the definition of $c_{\infty}$
and $b_{\lambda}$. ${\fim}$

The following result is the key point in the comparison of the
category of $\Omega$ with that of the sublevel sets of the
functional $J_{\lambda}$ given by $J_{\lambda}^{b_{\lambda,r}}=
\{u \in \mathcal{N}_{\lambda}:J_{\lambda}(u)\leq
b_{\lambda,r}\}$.

\begin{Le} \label{prop-baricentro}
There exists $\lambda_{*} > 0$ such that $\beta(u) \in
\Omega_{\lambda}^+$, whenever $u \in J_{\lambda}^{b_{r},\lambda}$
and $\lambda  \geq \lambda_{*}$.
\end{Le}

\noindent {\bf Proof.} The result follows from Lemma
\ref{lema-assintotico2} and arguments that were used in
\cite[Proposition 4.3]{Alv1}. ${\fim}$

Replacing $\Omega_{\lambda}$ by $B_{\lambda r}(0)$ we can prove that
there exists $u_{\lambda,r}$ a solution of problem
$$
 \left\{
\begin{array}{rcl}
\mathcal{L}u=f(u), \ B_{\lambda r}(0)\\
u>0, \ B_{\lambda r}(0)\\
u=0, \ \partial B_{\lambda r}(0)
\end{array}
\right.\leqno{(P_{\lambda,r})}
$$
which is radially symmetric on the origin (see \cite{BerestNirenb}).
For $\lambda>0$ and $r>0$, we define the operator
$\Phi_{r}:\lambda\Omega_{-}\rig H_{0}^{1}(\Omega_{\lambda})$ given
by
$$
 \Phi_{r,y}(x)=\left\{
\begin{array}{rcl}
u_{\lambda,r}(|x-y|), \ x\in B_{\lambda r}(y)\\
0, \ x\in \Omega_{\lambda}\backslash B_{\lambda r}(y)
\end{array}
\right.
$$

Notice that for each $y\in \lambda\Omega_{-}$, we have
$\beta(\Phi_{r,y})=y$.

\begin{Le}\label{75}
For each $\lambda\in [\lambda_{\ast},\infty)$, we have
$$
cat_{\Psi_{\lambda}^{b_{\lambda, r}}} \left( m_{\lambda}^{-1}(\Phi_{r}(\lambda\Omega_{-}))\right)=cat _{I_{\lambda}^{b_{\lambda, r}}} \left(\Phi_{r}(\lambda\Omega_{-})\right)\geq cat\Omega,
$$
where $\Psi_{\lambda}^{b_{\lambda, r}}=\{ u\in S_{\lambda}^{+}: \Psi_{\lambda}(u)\leq b_{\lambda, r}\}$.
\end{Le}
\noindent {\bf Proof.} Suppose that
$$
\Phi_{r}(\lambda\Omega_{-})=\bigcup_{j=1}^{k}A_{j}.
$$
where $A_{j}\subset I_{\lambda}^{b_{\lambda, r}}\subset \mathcal{N}_{\lambda}$ 
is closed and contractible in $I_{\lambda}^{b_{\lambda, r}}$. Since that
$m_{\lambda}:S_{\lambda}^{+}\rig \mathcal{N}_{\lambda}$ is a
homeomorphism we get
$$
cat_{\Psi_{\lambda}^{b_{\lambda, r}}} \left( m_{\lambda}^{-1}(\Phi_{r}(\lambda\Omega_{-}))\right)=cat _{I_{\lambda}^{b_{\lambda, r}}} \left(\Phi_{r}(\lambda\Omega_{-})\right)=k.
$$

Now the  rest of the proof follows from the Lemma \ref{prop-baricentro} and from \cite[Proposition
4.5]{Alv1}. $\fim$

\section{Proof of Theorem \ref{Theorem 1.1}}

Firstly, we define the compact set $K:= m_{\lambda}^{-1}(\Phi_{r}(\lambda\Omega_{-}))$ and we observe that 
$K\subset \Psi_{\lambda}^{b_{\lambda, r}}\subset S_{\lambda}^{+}$. Moreover, follow from the Lemma
$\ref{75}$ that 
$$
cat\Omega\leq cat_{\Psi_{\lambda}^{b_{\lambda, r}}}K.
$$
Follow from [\cite{Szulkin}, Theorem 27], with $c=b_{\lambda}<b_{\lambda, r}=d$, that $\Psi_{\lambda}^{b_{\lambda, r}}$
contains $cat \Omega$ critical points of $\Psi_{\lambda}$. From the Remark \ref{remark}, we conclude that $I_{\lambda}$
has at least $cat \Omega$ critical points, with energy in $[b_{\lambda}, b_{\lambda, r}]$.

On the other hand, if $cat\Omega>1$ we argue similarly to [\cite{BC 2}, Theorem 1.1]. Choosing $u^{\ast}\in S_{\lambda}^{+}$ such that $b_{\lambda, r}<\Psi_{\lambda}(u^{\ast})$,
we define
$$
\Theta=\{ t u^{\ast}+(1-t)u: t\in [0, 1] \ \mbox{e} \ u\in K\}.
$$
We observe that $\Theta$ is compact and $0\not\in \Theta$.

We also define the set
$$
\Gamma=\{\frac{w}{\|w\|_{\lambda}}: w\in \Theta\}\subset S_{\lambda}^{+}.
$$
Once $K\subset \Gamma$, with $K$ contractible in $\Gamma$ and
$$
b_{\lambda, r}<\Psi_{\lambda}(u^{\ast})\leq\max_{v\in \Gamma}\Psi_{\lambda}=:c,
$$
it follows that $K$ is contractible in $\Psi^{c}_{\lambda}$. Over again, from [\cite{Szulkin}, Theorem 27], with $2\leq k=cat \Omega$ e $e=c$,  we conclude that $\Psi_{\lambda}$ has
another critical point in $\Psi_{\lambda}^{c}\backslash \Psi_{\lambda}^{b_{\lambda, r}}$. From Remark \ref{remark},
it follows that $I_{\lambda}$ has another critical point with energy in $(b_{\lambda, r}, c]$.
$\fim$

\vspace{0.5cm} \noindent {\bf Acknowledgement.} This work was done
while the author was visiting the "Departamento de ecuaciones
diferenciales y an\'{a}lisis num\'{e}rico" of the Universidad de
Sevilla. They would like to express his gratitude to the Prof.
Antonio Suarez for his warm hospitality.

\end{document}